\documentclass[11pt]{article}

\usepackage{amsmath, amssymb, amscd}
\usepackage[matrix,arrow]{xy}

\def\eqref#1{(\ref{#1})}
\newcommand{\goth}{\frak}

\newcommand{\arrow}{{\:\longrightarrow\:}}
\newcommand{\Z}{{\Bbb Z}}
\newcommand{\C}{{\Bbb C}}
\newcommand{\R}{{\Bbb R}}

\renewcommand{\H}{{\Bbb H}}
\newcommand{\6}{\partial}
\def\1{\sqrt{-1}\:}
\newcommand{\restrict}[1]{{\left|_{{\phantom{|}\!\!}_{#1}}\right.}}

\renewcommand{\c}[1]{{\cal #1}}

\renewcommand{\tilde}{\widetilde}
\renewcommand{\bar}{\overline}
\renewcommand{\phi}{\varphi}
\renewcommand{\epsilon}{\varepsilon}
\renewcommand{\geq}{\geqslant}
\renewcommand{\leq}{\leqslant}

\newcommand{\even}{{\rm even}}

\newcommand{\End}{\operatorname{End}}
\newcommand{\Tot}{\operatorname{Tot}}

\renewcommand{\Re}{\operatorname{Re}}
\renewcommand{\Im}{\operatorname{Im}}

\newcommand{\comment}[1]{{}}

\def\blacksquare{\hbox{\vrule width 4pt height 4pt depth 0pt}}
\def\endproof{\blacksquare}

\makeatletter

\@ifundefined{Bbb}
     {\newcommand{\Bbb}[1]{{\mathbb #1}}}%
{}

\newcommand{\ps@verbit}{%
  \renewcommand{\@oddhead}{%
          \scriptsize
          {Examples HKT-manifolds}
          \hfil\tiny {M. Verbitsky, 11 March 2003 }}
  \renewcommand{\@evenhead}{\@oddhead}
  \renewcommand{\@oddfoot}{\hfil\thepage\hfil}
  \renewcommand{\@evenfoot}{\@oddfoot}}
 
\newcounter{Mycounter}[section]
\newcounter{lemma}[section]
\renewcommand{\thelemma}{{Lemma \thesection.\arabic{lemma}}}
\newcommand{\lemma}{%
     \setcounter{lemma}{\value{Mycounter}}
     \refstepcounter{lemma}
     \stepcounter{Mycounter}
     {\bf \thelemma:\ }}

\newcounter{claim}[section]
\renewcommand{\theclaim}{{Claim \thesection.\arabic{claim}}}
\newcommand{\claim}{%
     \setcounter{claim}{\value{Mycounter}}
     \refstepcounter{claim}
     \stepcounter{Mycounter}
     {\bf \theclaim:\ }}

\newcounter{sublemma}[section]
\newcounter{corollary}[section]
\renewcommand{\thecorollary}{{Corollary \thesection.\arabic{corollary}}}
\newcommand{\corollary}{%
     \setcounter{corollary}{\value{Mycounter}}
     \refstepcounter{corollary}
     \stepcounter{Mycounter}
     {\bf \thecorollary:\ }}

\newcounter{theorem}[section]
\renewcommand{\thetheorem}{{Theorem \thesection.\arabic{theorem}}}
\newcommand{\theorem}{%
     \setcounter{theorem}{\value{Mycounter}}
     \refstepcounter{theorem}
     \stepcounter{Mycounter}
     {\bf \thetheorem:\ }}

\newcounter{conjecture}[section]
\newcounter{proposition}[section]
\renewcommand{\theproposition}
       {{Proposition \thesection.\arabic{proposition}}}
\newcommand{\proposition}{%
     \setcounter{proposition}{\value{Mycounter}}
     \refstepcounter{proposition}
     \stepcounter{Mycounter}
     {\bf \theproposition:\ }}

\newcounter{definition}[section]
\renewcommand{\thedefinition}
       {{Definition~\thesection.\arabic{definition}}}
\newcommand{\definition}{%
     \setcounter{definition}{\value{Mycounter}}
     \refstepcounter{definition}
     \stepcounter{Mycounter}
     {\bf \thedefinition:\ }}

\newcounter{example}[section]
\renewcommand{\theexample}{{Example \thesection.\arabic{example}}}
\newcommand{\example}{%
     \setcounter{example}{\value{Mycounter}}
     \refstepcounter{example}
     \stepcounter{Mycounter}
     {\bf \theexample:\ }}

\newcounter{remark}[section]
\newcounter{problem}[section]
\newcounter{question}[section]
\@addtoreset{equation}{section}
\@addtoreset{footnote}{section}
\makeatother

\begin{document}

\begin{center}
{\LARGE\bf
Hyperk\"ahler manifolds with torsion \\[3mm]
obtained from hyperholomorphic bundles
}
\\[4mm]
Misha Verbitsky,\footnote{The author is 
partially supported by CRDF grant RM1-2354-MO02 and EPSRC grant  GR/R77773/01}
\\[4mm]
{\tt verbit@maths.gla.ac.uk, \ \  verbit@mccme.ru}
\end{center}

{\small 
\hspace{0.15\linewidth}
\begin{minipage}[t]{0.7\linewidth}
{\bf Abstract} \\
We construct examples of compact hyperk\"ahler
manifolds with torsion (HKT manifolds) which are
not homogeneous and not locally conformal hyperk\"ahler.
Consider a total space $T$ of a tangent bundle over
a hyperk\"ahler manifold $M$. The manifold $T$ is hypercomplex,
but it is never hyperk\"ahler, unless $M$ is flat.
We show that $T$ admits an HKT-structure. We also
prove that a quotient of $T$ by a $\Z$-action 
$v \arrow q^n v$ is HKT, for any real number $q\in \R$, 
$q>1$. This quotient is compact, if $M$ is compact.
A more general version of this construction holds 
for all hyperholomorphic bundles with holonomy 
in $Sp(n)$.
\end{minipage}
}

{
\small
\tableofcontents
}

\section{Introduction}
\label{_Intro_Section_}


Hyperk\"ahler manifolds with torsion (HKT-manifolds) were
introduced by P.S.Howe and G.Papadopoulos (\cite{_Howe_Papado_})
and much discussed in physics literature since then.
For an excellent survey of these works written from  a mathematician's
point of view, the reader is referred to the paper of
G. Grantcharov and Y. S. Poon \cite{_Gra_Poon_}.
In physics, HKT-manifolds appear as moduli of
brane solitons in supergravity and M-theory (\cite{_GP2_}, 
\cite{_P:Rome_}). HKT-manifolds also arise as moduli
space of some special black holes in N=2 supergravity 
(\cite{_GP1_}, \cite{_GPS_}).

The term ``hyperk\"ahler manifold with torsion'' is actually
quite misleading, because an HKT-manifold is not hyperk\"ahler.
This is why we prefer to use the abbreviation ``HKT-manifold''.

\hfill

HKT-manifolds are hypercomplex manifolds equipped with 
a special kind of Riemannian metrics.

A {\bf hypercomplex manifold} (\cite{_Boyer_})
is a $C^\infty$-manifold $M$ endowed
with a triple of almost complex structures $I, J, K\in \End(TM)$ 
which are integrable and satisfy the quaternionic relations
$I\circ J = - J\circ I = K$. If, in addition, 
$M$ is equipped with a Riemannian structure $g$
preserved by $I, J, K$, then $M$ is called
{\bf hypercomplex Hermitian}. If $(M, g)$ is
K\"ahler with respect to $I, J, K$, then
$(M, g, I, J, K)$ is called {\bf hyperk\"ahler}.

An HKT-manifold is a hypercomplex Hermitian
manifold which satisfies a similar, but weaker
condition \eqref{_HKT_intro_Equation_}. 

\hfill

Let $(M, g, I, J, K)$ be a hypercomplex Hermitian manifold.
Write the standard Hermitian forms on $M$ as follows:
\[ 
   \omega_I:= g(\cdot, I \cdot), \ \ \ \omega_J:= g(\cdot, J \cdot), \ \ \ 
   \omega_K:= g(\cdot, K \cdot).
\]
By definition, $M$ is hyperk\"ahler iff these forms are closed.
The HKT condition is weaker:
\begin{equation}\label{_HKT_intro_Equation_}
\6(\omega_J + \1 \omega_K)=0.
\end{equation}
Notice that $ \Omega=\frac 1 2(\omega_J + \1 \omega_K)$
is a $(2,0)$-form, for any hypercomplex Hermitian manifold,
as an elementary linear-algebraic calculation insures.
This form is called {\bf the canonical (2,0)-form associated with the
hypercomplex Hermitian structure}. As we shall see 
(\ref{_hh-metri_from_Omega_Proposition_}),
the metric can be recovered from the hypercomplex 
structure and the form $\Omega$.

Originally, the HKT-manifolds were defined in terms
of a quaternionic invariant connection with 
totally antisymmetric torsion (see \cite{_Howe_Papado_}, 
\cite{_Gra_Poon_}).

Many homogeneous examples of 
compact HKT-manifolds were obtained in 
\cite{_Howe_Papado_} and \cite{_Gra_Poon_}.
In \cite{_Ivanov:HKT_} it was shown that
any locally conformally hyperk\"ahler manifold
also admits an HKT-structure (see \cite{_Ornea:LCHK_}).

Locally, the HKT metrics can be studied using potential
functions (\cite{_Gra_Poon_}) in the same fashion
as one uses plurisubharmonic functions to study
K\"ahler metrics. This way one obtains many examples
of HKT-structures on a sufficiently small open
hypercomplex manifolds.

If $\dim_\R M=4$, every hypercomplex Hermitian metrics
is also HKT (the condition \eqref{_HKT_intro_Equation_}
is satisfied vacuously because the left hand side of
\eqref{_HKT_intro_Equation_} is a $(3,0)$-form). 

If $\dim_\R M>4$, the HKT-condition becomes highly non-trivial.
There are examples of hypercomplex manifolds not
admitting an HKT-structure (\cite{_Fino_Gra_}).
All known examples of compact HKT-manifolds are 
either homogeneous or locally conformally 
hyperk\"ahler.

\hfill

In the present paper, we construct
HKT-structures on fibered spaces associated with hyperk\"ahler
manifolds. A typical example of our construction is the following

\hfill

\theorem \label{_tot_HKT_Intro_Theorem_}
Let $M$ be a hyperk\"ahler manifold and 
\[ T^\circ M = \Tot(TM) \backslash \{\text{zero section}\}
\]
the total space of non-zero vectors in $TM$. Given $q \in \R$,
$|q|\neq 1$, let $\sim_q$ be the equivalence relation
generated by $x \sim_q qx$, $x \in TM$. Consider the 
quotient $T^\circ M/\sim_q$. Then $T^\circ M/\sim_q$
is equipped with a natural HKT-structure.

\hfill

{\bf Proof:} See \ref{_total_spa_Hopf_fibra_HKT_Theorem_}. \endproof

\hfill

\ref{_tot_HKT_Intro_Theorem_} is a special case of a 
much more general construction performed in 
Section \ref{_Hopf_bun_HKT_Section_}.


\section{The $q$-Dolbeault bicomplex}


In this Section, we 
introduce some notions of quaternionic linear algebra
which will be used further on. A reader well versed
in quaternions can safely skip this section.
We follow \cite{_V:projective_}.

\hfill

Let $M$ be a hypercomplex manifold, and 

\[ \Lambda^0 M \stackrel d\arrow \Lambda^1 M 
   \stackrel d\arrow \Lambda^2 M \stackrel d\arrow ...
\]
its de Rham complex. Consider the natural action of $SU(2)$ on 
$\Lambda^*M$.  Clearly, $SU(2)$ acts on 
$\Lambda^iM$, $i\leq \frac 1 2 \dim_\R M$
with weights $i, i-2, i-4, \dots$

We denote by $\Lambda^i_+$ the maximal $SU(2)$-subspace 
of $\Lambda^i$, on which $SU(2)$ acts with weight $i$. 

\hfill

The following linear-algebraic lemma allows one to compute
$\Lambda^i_+$ explicitly

\hfill

\lemma\label{_Lambda_+_explicit_Lemma_}
In the above assumptions, let $I$ be an induced complex structure,
and $\H_I$ the quaternion space, considered as a 2-dimensional
complex vector space with the complex structure induced by $I$.
Denote by $\Lambda^{p,0}_I(M)$ the space of $(p,0)$-forms on $(M,I)$. 
The space $\H_I$ is equipped with a
natural action of $SU(2)$. Consider $\Lambda^{p,0}_I(M)$
as a representation of $SU(2)$, with trivial group action. 
Then, there is a canonical isomorphism
\begin{equation}\label{_Lambda_+_explicit_Equation_}
\Lambda^p_+(M) \cong S^p_\C \H_I \otimes_\C \Lambda^{p,0}_I(M),
\end{equation}
where $S^p_\C \H_I$ denotes a $p$-th symmetric power of $\H_I$.
Moreover, the $SU(2)$-action on $\Lambda^p_+(M)$ is compatible with
the isomorphism \eqref{_Lambda_+_explicit_Equation_}.

\hfill

{\bf Proof:} This is \cite{_V:projective_}, Lemma 8.1.
\endproof

\hfill

Consider an $SU(2)$-invariant decomposition
\begin{equation}\label{_decompo_Lambda_to_Lambda_+_Equation_}
\Lambda^p(M) = \Lambda^p_+(M)\oplus V^p,
\end{equation}
where $V^p$ is the sum of all $SU(2)$-subspaces 
of $\Lambda^p(M)$ of weight less than $p$. Using the decomposition
\eqref{_decompo_Lambda_to_Lambda_+_Equation_},
we define the quaternionic Dolbeault differential 
$d_+:\; \Lambda^*_+(M)\arrow \Lambda^*_+(M)$ 
as a composition of de Rham differential and projection of to
$\Lambda^*_+(M) \subset \Lambda^*(M)$. Since the 
de Rham differential cannot increase the $SU(2)$-weight
of a form more than by 1, $d$ preserves
the subspace $V^*\subset \Lambda^*(M)$.
Therefore, $d_+$ is a differential in 
$\Lambda^*_+(M)$.

\hfill

Let $M$ be a hypercomplex manifold, and $I$ an induced complex
structure.  Consider the operator 
${\cal I}:\; \Lambda^*(M) \arrow \Lambda^*(M)$
mapping a $(p,q)$-form $\eta$ to $\1(p-q)\eta$.
By definition, ${\cal I}$ belongs to the Lie algebra $\goth{su}(2)$
acting on $\Lambda^*(M)$ in the standard way. 
Therefore, ${\cal I}$ preserves the subspace
$\Lambda^*_+(M) \subset \Lambda^*(M)$. 
We obtain the Hodge decomposition
\[ \Lambda^*_+(M) = \oplus_{p,q} \Lambda^{p,q}_{+,I}(M). \]

\hfill

Let $M$ be a hypercomplex manifold, $I$ an induced comlex structure,
and $I, J, K\in \H$ the standard triple of induced complex structures. 
Clearly, $J$ acts on the complexified co tangent space
$\Lambda^1M\otimes \C$ mapping $\Lambda_I^{0,1}(M)$
to $\Lambda_I^{1,0}(M)$. 
Consider a differential operator 
\[ \6_J:\; C^\infty(M)\arrow \Lambda_I^{1,0}(M),\]
mapping $f$ to $J(\bar\6 f)$, where 
$\bar\6:\; C^\infty(M)\arrow \Lambda_I^{0,1}(M)$
is the standard Dolbeault differential on a K\"ahler
manifold $(M,I)$. We extend $\6_J$
to a differential
\[ 
   \6_J:\; \Lambda_I^{p,0}(M)\arrow \Lambda_I^{p+1,0}(M),
\]
using the Leibniz rule.

\hfill

\proposition \label{_d_+_Hodge_6_J_Proposition_}
Let $M$ be a hypercomplex manifold, $I$ an induced 
complex structure, $I, J, K$ the standard 
basis in quaternion algebra, and
\[ \Lambda^*_+(M)= \oplus_{p,q} \Lambda^{p,q}_{I,+}(M) \]
the Hodge decomposition of the quaternionic Dolbeault
complex. Then there exists a canonical isomorphism
\begin{equation}\label{_Lambda_+_Hodge_deco_expli_Equation_} 
   \Lambda^{p,q}_{I, +}(M)\cong \Lambda^{p+q, 0}_I(M). 
\end{equation}
Under this identification, the quaternionic Dolbeault differential
\[ d_+:\; \Lambda^{p,q}_{I, +}(M)\arrow 
   \Lambda^{p+1,q}_{I, +}(M)\oplus \Lambda^{p,q+1}_{I, +}(M)
\]
corresponds to a sum
\[ \6 \oplus \6_J:\; \Lambda^{p+q, 0}_{I}(M) \arrow
   \Lambda^{p+q+1, 0}_{I}(M)\oplus \Lambda^{p+q+1,0}_{I}(M).
\]

{\bf Proof:} This is Proposition 8.13 of \cite{_V:projective_}.
\endproof

\hfill

The statement of 
\ref{_d_+_Hodge_6_J_Proposition_}
can be represented 
by the following diagram

\begin{equation}\label{_bicomple_XY_Equation}
\begin{minipage}[m]{0.85\linewidth}
{\tiny $
\xymatrix @C+1mm @R+10mm@!0  { 
  {} &{} & \Lambda^0_+(M) \ar[dl]^{d'_+} \ar[dr]^{d''_+} 
   {} &{} & {} &{} & {} &{} & \Lambda^{0,0}_I(M) \ar[dl]^{ \6} \ar[dr]^{ \6_J}
   {} &{} &  \\
 {} & \Lambda^{1,0}_+(M) \ar[dl]^{d'_+} \ar[dr]^{d''_+} {} &
 {} & \Lambda^{0,1}_+(M) \ar[dl]^{d'_+} \ar[dr]^{d''_+}{} &{} & 
\text{\large $\cong$} {} &
 {} &\Lambda^{1,0}_I(M)\ar[dl]^{ \6} \ar[dr]^{ \6_J}{} &  {} &
 \Lambda^{1,0}_I(M)\ar[dl]^{ \6} \ar[dr]^{ \6_J}{} &\\
 \Lambda^{2,0}_+(M) {} &{} & \Lambda^{1,1}_+(M) 
   {} &{} & \Lambda^{0,2}_+(M){} & \ \ \ \ \ \ {} & \Lambda^{2,0}_I(M){} & {} & 
\Lambda^{2,0}_I(M) {} & {} & \Lambda^{2,0}_I(M) \\
}
$
}
\end{minipage}
\end{equation}
where $d_+= d'_+ + d''_+$ is the Hodge decomposition of 
the quaternionic Dolbeault differential.

\hfill

Using the $SU(2)$-action, we may identify the
bundles $\Lambda^{p,q}_+(M)$ with 
$\Lambda^{p+q,0}_+(M) = \Lambda^{p+q,0}_I(M)$
explicitly, as follows.

Let ${\cal J}, {\cal K}$ be the Lie algebra operators acting on 
differential forms and associated
with $J$, $K$ 
in the same way as ${\cal I}$ is associated with $I$.
Consider the map ${\cal R}:\; \Lambda^*(M) \arrow \Lambda^*(M)$,
\begin{equation} \label{_c_d_Definition_Equation_}
  {\cal R} := \frac{{\cal J} - \1 {\cal K}}{2}.
\end{equation}
It is easy to check that the Lie algebra elements
${\cal R}, {\cal I}, \bar{{\cal R}}$
form an $SL(2)$-triple in the complexification
of the standard $SU(2)\subset \End(\Lambda^*(M))$.
Therefore, ${\cal R}$ maps $\Lambda^{p,q}_+(M)$
to $\Lambda^{p+1,q-1}_+(M)$.
Since $\Lambda^m_+(M)$ is a representation of
weight $m$, ${\cal R}$ induces an isomorphism
\[ {\cal R}:\; \Lambda^{p,q}_{+,I}(M)\arrow \Lambda^{p+1,q-1}_{+,I}(M),
\]
for all $q>0$.

Together with 
\eqref{_bicomple_XY_Equation},
this observation implies the following.

\hfill

\claim\label{_c_R_on_2-forms_Claim_}
Let $M$ be a hypercomplex manifold,
$I$ an induced complex structure,  
and $\eta\in \Lambda^{1,1}_I(M)$ a (1,1)-form. 
Then $\eta$ is $SU(2)$-invariant if and only if
${\cal R}(\eta)=0$. Moreover, for all functions
$\psi$ on $M$, we have
\[ 
{\cal R}(\6\bar\6(\psi)) = \6\6_J(\psi).
\]

\endproof

\hfill

Assume now that the manifold $M$ is hypercomplex Hermitian.
Consider the 3-dimensional space generated by the 2-forms
$\omega_I$, $\omega_J$ and $\omega_K$. This is a weight 2
representationn of $SU(2)$. Moreover,
that 
\begin{equation} \label{_R(omega)_Equation_}
 {\cal R}(\omega_I) = \Omega,
\end{equation}
where
$\Omega=\frac 1 2(\omega_J + \1 \omega_K)$ 
is the canonical $(2,0)$-form.


\section{The $q$-positive forms}


Let $M$ be a hypercomplex manifold, and
$\Lambda^{p,q}_I(M)$ the bundle of $(p,q)$-forms
on $(M,I)$. Consider the map
$J:\; \Lambda^*(M)\arrow \Lambda^*(M)$,
\[ 
J(dx_1\wedge dx_2 \wedge ...) = J(dx_1) \wedge J(dx_2) \wedge ...
\]
Clearly, on 2-forms we have $J^2 =1$; more generally,
\begin{equation}\label{_J^2=1_Equation_}
\left(J\restrict{\Lambda^\even(M)}\right)^2 =1.
\end{equation}
Since $J$ and $I$ anticommute, we have
$J(\Lambda^{p,q}_I(M)) = \Lambda^{q,p}_I(M)$.
By \eqref{_J^2=1_Equation_}, the map
$\eta\arrow J(\bar \eta)$ defines a real structure
on $\Lambda^{2,0}_I(M)$.

\hfill

\definition
Let $\eta\in \Lambda^{2,0}_I(M)$ be a (2,0)-form
on a hypercomplex manifold $M$.
We say that $\eta$ is $q$-real if $\eta = J(\bar\eta)$.
We say that $\eta$ is $q$-positive if $\eta$ is $q$-real, and
\begin{equation}\label{_q-posi_Equation_}
\eta(v, J(\bar v))\geq 0
\end{equation}
for any $v \in T^{1,0}_I(M)$.
We say that $\eta$ is {\bf strictly $q$-positive}
if the inequality \eqref{_q-posi_Equation_}
is strict, for all $v\neq 0$.

The $q$-positive forms were introduced
and studied at some length in \cite{_Verbitsky:Hyperholo_stable_},
under the name ``$K$-positive forms''.
These forms were used to study the stability 
of certain coherent sheaves. Some properties of
$q$-positive forms are remarkably close to that
of the usual positive forms, studied in algebraic
geometry in connection with Vanishing Theorems.

\hfill

\proposition\label{_hh-metri_from_Omega_Proposition_}
Let $M$ be a hypercomplex manifold, and $h$  
a hypercomplex Hermitian metric. Consider the form
\[ \Omega:= \omega_J + \1 \omega_K\]
(see \eqref{_HKT_intro_Equation_}).
Then $\Omega$ is strictly $q$-positive.
Conversely, every strictly $q$-positive
$(2,0)$-form is obtained from a unique
hypercomplex Hermitian metric on $M$.

\hfill

{\bf Proof:} The form $\Omega:= \omega_J + \1 \omega_K$
is $q$-positive as an elementary calculation insures.
Indeed, write the orthonormal basis
$\xi_1, \xi_2, ... \xi_{2n} \in \Lambda^{1,0}(M)$
in such a way that
\begin{equation}\label{_J_acts_on_xi_Equation_}
J(\xi_{2i-1})= \bar\xi_{2i},\ \ \  J(\xi_{2i})= -\bar\xi_{2i-1}. 
\end{equation}
Then 
\begin{equation}\label{_Omega_expli_Equation_}
\Omega= \xi_1\wedge \xi_2 + \xi_3\wedge \xi_4 + ...
\end{equation}
This form is clearly $q$-real and strictly $q$-positive.

Conversely, let $\Omega$ be a $q$-real and strictly $q$-positive
form on a hypercomplex manifold $M$. We can write $\Omega$
is coordinates as
\[
\Omega= \alpha_1\xi_1\wedge \xi_2 + \alpha_3\xi_3\wedge \xi_4 + ...
\] 
where $\alpha_i$ are positive real numbers, and $\xi_i$ satisfy
\eqref{_J_acts_on_xi_Equation_}. 

Write a hypercomplex Hermitian form $h$ as
\begin{equation}\label{_hc-H_via_Omega_Equation_}
\begin{aligned}
 h = & \alpha_1 ((\Re \xi_1)^2 + (\Im \xi_1)^2+ (\Re \xi_2)^2 + (\Im \xi_2)^2)\\
 {\ }  &+ \alpha_3 ((\Re \xi_3)^2 + (\Im \xi_3)^2+ (\Re \xi_4)^2 + (\Im \xi_4)^2) + ...
\end{aligned}
\end{equation}
Clearly, the corresponding canonical $(2,0)$-form is 
equal $\Omega$.

The Hermitian metric \eqref{_hc-H_via_Omega_Equation_} 
can be reconstructed from $\Omega$ directly as follows:
\[
h(x, y) = \Omega(x^{1,0}, J(y^{0,1})),
\]
for all $x, y \in T_\R M$,
where $x^{1,0}$, $y^{0,1}$ denotes the
$(1,0)$ and $(0,1)$-parts of $x$, $y$.
We proved that the hypercomplex Hermitian structure
is uniquely determined by the strictly $q$-positive
form $\Omega$. \endproof

\hfill

The following Corollary gives an interpretation
of HKT-structures in terms of the canonical $(2,0)$-form. 

\hfill

\corollary \label{_HKT_via_q-pos_Corollary_}
Let $M$ be a hypercomplex manifold,
and $\Omega\in \Lambda^{2,0}_(M)$
a $\6$-closed strictly $q$-positive $(2,0)$-form. Then
$M$ is an HKT-manifold, and $\Omega$ is obtained
as a canonical $(2,0)$-form of an HKT-metric $h$.
Moreover, $h$ is uniquely determined by $\Omega$.

\hfill

{\bf Proof:}  By \ref{_hh-metri_from_Omega_Proposition_},
$\Omega= \omega_J + \1 \Omega_K$, for some
hypercomplex Hermitian metric $h$. 
Since $\6 \Omega=0$, $(M, h)$ is an HKT-manifold.
\endproof


\section{Hyperholomorphic bundles}


Hyperholomorphic bundles were introduced and studied at some
length in \cite{_Verbitsky:Hyperholo_bundles_}.
 Let $B$ be a holomorphic vector bundle over a complex
manifold $X$, $\nabla$ a  connection 
in $B$ and $\Theta\in\Lambda^2\otimes End(B)$ be its curvature. 
This connection
is called {\bf compatible with the holomorphic structure} if
$\nabla_\gamma(\zeta)=0$ for any holomorphic section $\zeta$ and
any antiholomorphic tangent vector field $\gamma\in T^{0,1}(X)$. 
If there exists a holomorphic structure compatible with the given
Hermitian connection then this connection is called
{\bf integrable}.

\hfill

\theorem \label{_Newle_Nie_for_bu_Theorem_}
Let $\nabla$ be a Hermitian connection in a complex vector
bundle $B$ over a complex manifold $X$. Then $\nabla$ is integrable
if and only if $\Theta\in\Lambda^{1,1}(X, \End(B))$, where
$\Lambda^{1,1}(X, \End(B))$ denotes the forms of Hodge
type (1,1). Also,
the holomorphic structure compatible with $\nabla$ is unique.

{\bf Proof:} This is Proposition 4.17 of \cite{_Kobayashi_}, 
Chapter I.
\endproof

\hfill

This proposition is a version of Newlander-Nirenberg theorem.
For vector bundles, it was proven by M. Atiyah and R. Bott.

\hfill

\definition \label{_hyperho_conne_Definition_}
Let $B$ be a Hermitian vector bundle with
a connection $\nabla$ over a hypercomplex manifold
$M$. Then $\nabla$ is called {\bf hyperholomorphic} if 
$\nabla$ is
integrable with respect to each of the complex structures induced
by the hypercomplex structure. 
 
As follows from 
\ref{_Newle_Nie_for_bu_Theorem_}, 
$\nabla$ is hyperholomorphic
if and only if its curvature $\Theta$ is of Hodge type (1,1) with
respect to any of the complex structures induced by a hypercomplex 
structure.

An easy calculation shows that
$\nabla$ is hyperholomorphic
if and only if $\Theta$ is an 
$SU(2)$-invariant differential form.

\hfill

Hyperholomorphic 
bundles are quite ubiquitous.
Clearly, the tangent bundle to a hyperk\"ahler manifold 
and all its tensor powers are hyperholomorphic. 
There are many other examples

\hfill

\example
Let $M$ be a compact hyperk\"ahler manifold,
$B$ a holomorphic bundle. Then $B$ admits a
unique hyperholomorphic connection, if
$B$ is stable and the cohomology classes 
$c_1(B)$ and $c_2(B)$ are $SU(2)$-invariant.
Moreover, if $M$ is generic in its deformation
class, then all stable bundles admit
a hyperholomorphic connection.


\section{${\Bbb H}$-hyperholomorphic bundles}


\definition
Let $M$ be a hypercomplex manifold, and $(B, \nabla)$ a hyperholomorphic
bundle on $M$, $\dim_\C B = 2n$. The bundle
$B$ is called {\bf ${\Bbb H}$-hyperholomorphic}
if $\nabla$ preserves a $\C$-linear symplectic structure
on $B$. In other words, $B$ is {\bf ${\Bbb H}$-hyperholomorphic}
if the holonomy of $\nabla$ is contained in $Sp(n)$.

\hfill

The following examples are obvious.

\hfill

\example \label{_F+F^*_H-hh_Example_}
Let $F$ be a hyperholomorphic
bundle on $M$. Then $F \oplus F^*$ is ${\Bbb H}$-hyperholomorphic.

\hfill

\example \label{_TM_H-hh_Example_}
Consider the tangent bundle $TM$ on $M$. Assume that $M$
is hyperk\"ahler. Then $TM$ is
${\Bbb H}$-hyperholomorphic.

\hfill

The main property of ${\Bbb H}$-hyperholomorphic
bundles is the following. 

\hfill

\claim \label{_Tot_B_hc_Claim_}
Let $M$ be a hyperk\"ahler manifold, 
and $B$ an ${\Bbb H}$-hy\-per\-ho\-lo\-mor\-phic
bundle. Denote by $\Tot B$ the total space of $B$. Then 
$\Tot B$ is equipped with a natural hypercomplex structure.
In particular, the total space of $TM$ is hypercomplex. 

\hfill

{\bf Proof:} Since the holonomy of $B$ is contained in 
$Sp(n)$, there is a natural parallel action of $\Bbb H$ on $B$.
Given a quaternion $L\in \Bbb H$, $L^2=-1$, consider
$B$ as a complex vector bundle with the complex structure
defined by $L$. Denote this complex vector bundle as $(B,L)$.
Since the curvature of $B$
is $SU(2)$-invariant, the bundle $(B,L)$
is hyperholomorphic. Consider $(B,L)$
as a holomorphic vector bundle on $(M,L)$.
Denote the corresponding complex structure
on $\Tot B$ by $L$. We obtained an
integrable complex structure on $\Tot B$ 
for each quaternion $L\in \Bbb H$, $L^2=-1$.
It is easy to check that these
complex structures satisfy quaternionic
relations, inducing a hypercomplex
structure on $\Tot B$.
\endproof

\section{The Obata connection on $\Tot B$.}

Let $M$ be a hyperk\"ahler manifold, and $B$ an ${\Bbb H}$-hyperholomorphic
bundle. By \ref{_Tot_B_hc_Claim_}, the total space $\Tot B$
is hypercomplex. One can ask whether this hypercomplex structure
is hyperk\"ahler. The answer is - never (unless $B$ is flat).

\hfill

Given a hypercomplex manifold, one can easily 
establish whether $M$ admits a hyperk\"ahler structure.
This is done most easily using the so-called Obata connection.

\hfill

\theorem
(Obata)
Let $M$ be a hypercomplex manifold. Then $M$ admits
a unique torsion-free connection  which preserves the
hypercomplex structure.\footnote{This connection is called 
the Obata connection.}

\hfill

{\bf Proof:} Well known (see \cite{_Obata_}). \endproof

\hfill

If $M$ is hyperk\"ahler, then the Levi-Civita connection 
preserves the hypercomplex structure. In this case,
the Levi-Civita connection coincides with the Obata connection.

\hfill

To determine whether a hypercomplex manifold $M$ admits a 
hyperk\"ahler structure, one needs to compute the holonomy of
the Obata connection. The manifold is hyperk\"ahler
if and only if the holonomy $Hol$ preserves a metric;
that is, $M$ is hyperk\"ahler if and only if 
$Hol$ is contained in $Sp(n)$.

\hfill

\proposition
Let $M$ be a hyperk\"ahler manifold, $B$ an ${\Bbb H}$-hy\-per\-ho\-lo\-mor\-phic
bundle, and $\Tot B$ its total space considered
as a hypercomplex manifold (see 
\ref{_Tot_B_hc_Claim_}). Assume that the
curvature of $B$ is non-zero. Then
$\Tot B$ does not admit a hyperk\"ahler structure.

\hfill

{\bf Proof:} One could compute the holonomy group of the Obata connection
of $\Tot B$, and show that it is non-compact. To avoid 
excessive computations, we use a less straightforward 
argument. 

Suppose that $\Tot B$ is hyperk\"ahler. Given
$m\in M$, let $B_m\subset \Tot B$ be the fiber of $B$ in $m$.
By construction, $B_m$ is a hypercomplex submanifold in $\Tot B$.
Such submanifolds are called {\bf trianalytic} (see 
\cite{_Verbitsky:Symplectic_II_}, \cite{_Verbitsky:Deforma_}
for a study of trianalytic cycles on hyperk\"ahler manifolds).
In \cite{_Verbitsky:Deforma_}, it was shown that trianalytic
submanifolds are completely geodesic. In other words, 
for any trianalytic submanifold $Z \subset X$,
the Levi-Civita connection on $TX\restrict Z$
preserves the orthogonal decomposition
\begin{equation}\label{_ortho_decompo_tangent_Equation_} 
  TX\restrict Z = TZ \oplus TZ^\bot 
\end{equation}
If we have a hypercomplex fibration $X\arrow Y$,
the decomposition \eqref{_ortho_decompo_tangent_Equation_} 
gives a connection for this fibration. In \cite{_Verbitsky:Deforma_}
it was shown that this connection is flat, for any hyperk\"ahler
fibration. 

We obtain a flat connection $\nabla$ in the fibration
$\Tot B \arrow M$. This connection is clearly compatible
with the additive structure on the bundle $B$. Therefore,
$\nabla$ is an affine connection on $B$. By construction,
$\nabla$ is compatible with the hypercomplex structure on 
$\Tot B$. Therefore, $\nabla$ coincides with the 
hyperholomorphic connection on $B$. We proved
that $B$ is flat. \endproof


\section{HKT-structure on $\Tot B$.}
\label{_Tot_B_HKT_Section_}


Let $M$ be a smooth manifold. Given a bundle with connection
on $M$, we have a decomposition

\begin{equation} \label{_TTotB_Equation_}
  T \Tot B = T_{ver} \oplus T_{hor}
\end{equation}
of the tangent space to $\Tot B$ into horizontal and vertical components.
Clearly, the bundle $T_{ver}$ is identified with $\pi^* B$,
and $T_{hor}$ with $\pi^* TM$,
where $\pi:\; \Tot B \arrow M$ is the standard projection.

Assume now that $M$ is a Riemannian manifold, and $B$ 
a vector bundle, equipped with a Euclidean metric. Then $\Tot B$ is equipped
with a Riemannian metric $g$ defined by the following conditions.

\hfill

\begin{description}
\item[(i)] The decomposition 
$T \Tot B = T_{ver} \oplus T_{hor}$ is orthogonal with respect to $g$.

\item[(ii)] Under the natural identification $T_{ver}\cong \pi^* B$,
the metric $g$ restricted to $T_{ver}$ becomes the metric on $B$.

\item[(iii)] The metric $g$ restricted $T_{hor}\cong  \pi^*TM$ 
is equal to the metric induced on $\pi^* TM$ from the Riemannian structure
on $M$. 

\end{description}

\hfill

\definition \label{_sta_metric_on_total_spa_Definition_}
In the above assumptions, the metric $g$ is called {\bf the natural
metric on $\Tot B$ induced by the connection and the metrics 
on $M$ and $B$}.

\hfill

Notice that the metric $g$ depends from the metrics on $B$ and $M$ and 
from the connection in $B$. Different connections induce different
metrics on $\Tot B$.

\hfill

\theorem \label{_Tot_B_HKT_Theorem_}
Let $M$ be an HKT-manifold, and $B$ an ${\Bbb H}$-hy\-per\-ho\-lo\-mor\-phic
vector bundle on $M$. Consider the metric $g$ on $\Tot B$
defined as in \ref{_sta_metric_on_total_spa_Definition_}.
Then $g$ is an HKT-metric.

\hfill

{\bf Proof:} Consider the decomposition
$g = \pi^* g_M + \pi^* g_B$ of the metric $g$ 
onto the horizontal and vertical components. 
Since the decomposition
$T \Tot B = T_{ver} \oplus T_{hor}$
is compatible with the hypercomplex structure,
the 2-forms $g_{hor}:=\pi^* g_M$ and $g_{ver}:=\pi^* g_B$
are $SU(2)$-invariant. Consider the corresponding
(2,0)-forms $\Omega_{hor}$ and $\Omega_{ver}$
obtained as in \eqref{_HKT_intro_Equation_};
\[ 
  \Omega_{hor} = {\omega_J}_{hor} + \1{\omega_K}_{hor} 
\]
where ${\omega_J}_{hor} = g_{hor}(J\cdot, \cdot)$, 
${\omega_K}_{hor}= g_{hor}(K\cdot, \cdot)$ are
differential forms associated with $g_{hor}$ and $J, K$ as in
\eqref{_HKT_intro_Equation_}.

Then $\Omega_{hor}$ and $\Omega_{ver}$
are horizontal and vertical components of the 
standard (2,0)-form of $\Tot B$:
\begin{equation}\label{Omega_decompo_Equation_}
   \Omega= \Omega_{hor} +\Omega_{ver} 
\end{equation}
The HKT condition can be written as
$\6\Omega=0$ \eqref{_HKT_intro_Equation_}. 
Let $\Omega_M$ be the standard
$(2,0)$-form of $M$. Since $M$ is an HKT manifold, 
\eqref{_HKT_intro_Equation_} holds on $M$ and 
the form $\Omega_{hor}$
satisfies 
\[ 
    \6 \Omega_{hor} =  \6 \pi^*\Omega_M= 0.
\]
To prove \ref{_Tot_B_HKT_Theorem_},
it remains to show 
\begin{equation} \label{_6_Omega_ver_Equation_}
\6\Omega_{ver}=0
\end{equation}
Consider a function 
\begin{equation}\label{_length_Psi_Equation_} 
   \Psi:\; \Tot B \arrow \R, \ \ \ \Psi(v) = |v|^2, 
\end{equation}
mapping a vector $v\in TM$ to the square of 
its norm.  Let
\[
0 \arrow \Omega^{1,0} \stackrel{\6, \6_J} \arrow
  \Omega^{2,0} \stackrel{\6, \6_J} \arrow 
  \Omega^{3,0} \stackrel{\6, \6_J} \arrow ...
\]
be the bicomplex defined in \eqref{_bicomple_XY_Equation}.
To prove \eqref{_6_Omega_ver_Equation_}, 
and hence \ref{_Tot_B_HKT_Theorem_},
it suffices to prove
\begin{equation}\label{_pote_for_Omega_ver_Equation_}
\6\6_J \Psi= \Omega_{ver}.
\end{equation}
By \ref{_c_R_on_2-forms_Claim_}, we have
\[ \6\6_J \Psi  = {\cal R}(\6\bar\6 \Psi), 
\]
where $R:\; \Lambda^{1,1}(\Tot B) \arrow \Lambda^{2,0}(\Tot B)$
is the operator
\[ {\cal R} = \frac{{\cal J} - \1 {\cal K}}{2}
\]
(see \eqref{_c_d_Definition_Equation_}). 
However, the 2-form $\6\bar\6 \Psi$
is quite easy to compute. From \cite{_Besse:Einst_Manifo_},
(15.19), we obtain:
\begin{equation}\label{_6_bar_6_Psi_Equation_}
\6\bar\6 \Psi = \omega_{ver} + \xi,
\end{equation}
where $\omega_{ver}= g_{ver}(\cdot, I \cdot)$ 
is the Hermitian form of $g_{ver}$, and $\xi$ is defined as following.
Using the decomposition \eqref{_TTotB_Equation_}, 
we consider $\Lambda^2 T_{hor}$ as a subbundle in $\Lambda^2 \Tot B$.
Then  $\xi \in \Lambda^2 T_{hor} \subset \Lambda^2 \Tot B$
is a 2-form on $T_{hor}$ mapping a pair of vectors $(x, y)$
\begin{align*}
    x, y &\in T_{hor}\restrict{(m,b)} \subset T_{(m,b)}\Tot B, \\
   & T_{hor}\restrict{(m,b)} = T_m M, \\
   & (m,b)\in \Tot B, m\in M, b \in B\restrict m
\end{align*}
to $(R(x,y,b) \bar b)$, where $R\in \Lambda^2 M \otimes \End B$
is the curvature of $B$. The form $\xi$ is $SU(2)$-invariant
because the curvature of $B$ is $SU(2)$-invariant. Therefore,
${\cal R}(\xi)=0$ (\ref{_c_R_on_2-forms_Claim_}), and 
\begin{equation}\label{_6_6_J_Psi_Equation_}
  \6\6_J \Psi = {\cal R}(\6\bar\6 \Psi) =   
  {\cal R}(\omega_{ver}) = \Omega_{ver}
\end{equation}
(the last equation holds by \eqref{_R(omega)_Equation_}). This proves
\eqref{_pote_for_Omega_ver_Equation_}.
\ref{_Tot_B_HKT_Theorem_} is proven.
\endproof


\section{New examples of compact HKT-manifolds}
\label{_Hopf_bun_HKT_Section_}


Let $M$ be a compact HKT-manifold, e.g. a hyperk\"ahler
manifold, and $B$ an ${\Bbb H}$-hyperholomorphic vector bundle on $M$
(for examples of ${\Bbb H}$-hy\-per\-ho\-lo\-mor\-phic vector bundles
see \ref{_F+F^*_H-hh_Example_} and \ref{_TM_H-hh_Example_}).
Denote by $\Tot^\circ B$ be the space of non-zero vectors in $B$.
Fix a real number $q>1$. Consider the map 
\[ 
  \rho_q:\; \Tot^\circ B\arrow \Tot^\circ B, \ \ 
  \rho_q(b) = qb,\ \  b \in \Tot^\circ B,
\] 
and let $\c M = \Tot^\circ B/\rho_q$
be the corresponding quotient space. Since the map $b \arrow qb$
is compatible with the hypercomplex structure, the space
$\c M$ is hypercomplex.  It is fibered over
a compact manifold $M$, with fibers Hopf manifolds which are
homeomorphic to $S^1 \times S^{2m-1}$, $m = \dim_\R B$,
hence it is compact.

\hfill

\theorem \label{_total_spa_Hopf_fibra_HKT_Theorem_}
In the above assumptions, $\c M$ admits a natural HKT-structure.

\hfill

{\bf Proof:} By \ref{_HKT_via_q-pos_Corollary_}, 
we need to construct a $q$-positive $\6$-closed
$(2,0)$-form on $\c M$. Let $\tilde \Omega$ be a (2,0)-form
on $\Tot^\circ B$, 
\[ \tilde \Omega = \pi^* \Omega_M + \6\6_J \log \Psi,
\]
where $\pi^* \Omega_M$ is the canonical (2,0)-form
on $M$ lifted to $\Tot^\circ B$, and $\Psi:\; \Tot B \arrow \R$
the square norm function \eqref{_length_Psi_Equation_}.
The map $v \stackrel {\rho_q}\arrow qv$ satisfies
$\rho_q^* \log \Psi = \log \Psi + \log q^2$,
and therefore 
\[
\rho_q^*\6\6_J  \log \Psi = \6\6_J  \log \Psi.
\]
This implies that $\tilde \Omega =  \pi^* \Omega_M + \6\6_J \log \Psi$
is $\rho_q$-invariant, hence defines a form $\Omega$ on 
 $\c M = \Tot^\circ B/\rho_q$.

By construction, the form $\Omega$ is $\6$-closed. To prove
\ref{_total_spa_Hopf_fibra_HKT_Theorem_}, it remains to show 
that $\tilde \Omega$ is strictly $q$-positive. We use
the same argument as used to show that a locally conformal
hyperk\"ahler manifold is HKT. 

We have 
\begin{equation}\label{_6_6_J_log_Psi_Equation_}
\6\6_J \log \Psi = \frac{\6\6_J \Psi}{\Psi} - 
       \frac{\6\Psi\wedge \6_J \Psi}{\Psi^2}.
\end{equation}
In all directions orthogonal to $\6\Psi,\6_J \Psi$,
the form $\6\6_J \log \Psi$ is proportional to 
$\6\6_J \Psi$, hence $q$-positive by \eqref{_6_6_J_Psi_Equation_}.
Moreover, \eqref{_6_6_J_log_Psi_Equation_} implies that
\[
\tilde \Omega = \Omega_{hor} + \frac{\Omega_{ver}}{\Psi}
 - \frac{\6\Psi\wedge \6_J \Psi}{\Psi^2},
\] 
(we use the notation introduced in 
Section \ref{_Tot_B_HKT_Section_}).
The form $\Omega_{hor} + \frac{\Omega_{ver}}{\Psi}$
is strictly $q$-positive (\ref{_Tot_B_HKT_Theorem_}).
The vertical and the horizontal tangent vectors
are orthogonal with respect to $\tilde \Omega$.
Since $\frac{\6\Psi\wedge \6_J \Psi}{\Psi^2}$ vanishes on
all horizontal tangent vectors, it remains to prove that 
$\tilde \Omega(x, J \bar x)>0$, 
where $x$ is vertical.

Let $\xi\in T^{1,0}(\Tot^\circ B)$ be the vertical
tangent vector to $\Tot^\circ B$ which is dual to
$\frac{d \Psi}{\sqrt \Psi}$. Clearly, $\6 \Psi$ is the $(1,0)$-part of 
$\xi$. For all $x\in T^{1,0}(\Tot^\circ B)$, we have
\[
 \frac{\6\Psi\wedge \6_J \Psi}{\Psi^2}\bigg(x, J(\bar x)\bigg) = 
 (\xi, x)_{H}^2,
\]
where $(\cdot, \cdot)_{H}$ denotes the Riemannian form.
Similarly, 
\[ \Omega_{ver} (x, J(\bar x)) = 2(x,x)_{H}
\]
(this can be checked by writing $\Omega_{ver}$
is coordinates as in \eqref{_Omega_expli_Equation_}).
Using Cauchy inequality and $|\xi|=1$, we obtain
$(x,x)_{H} \geq (\xi, x)_{H}^2$.
Then 
\begin{align*}
\tilde \Omega(x, J(\bar x)) &= \frac{\Omega_{ver}}{\Psi}\bigg(x, J(\bar x)\bigg)
 - \frac{\6\Psi\wedge \6_J \Psi}{\Psi^2}\bigg(x, J(\bar x)\bigg)\\
 &= 2 \frac{(x,x)_{H}}\Psi - \frac{(\xi, x)_{H}^2}\Psi 
 \geq \frac{(x,x)_{H}}\Psi >0
\end{align*}
for all vertical tangent vectors $x\neq 0$.
This proves \ref{_total_spa_Hopf_fibra_HKT_Theorem_}.
\endproof

\hfill

{\bf Acknowledgements:}
I am grateful to G. Grantcharov,  D. Kaledin
and Y.-S. Poon for interesting discussions.

\hfill


\end{document}